\documentclass[12pt,a4paper]{article}
\usepackage[utf8]{inputenc}
\usepackage[T1]{fontenc}
\usepackage[english]{babel}
\usepackage{amsmath,amssymb,amsthm,booktabs}
\usepackage[margin=2.5cm]{geometry}
\usepackage{tikz}
\usetikzlibrary{arrows.meta,decorations.pathreplacing,matrix,calc}

\usepackage[colorlinks=true,linkcolor=blue,citecolor=blue,urlcolor=blue]{hyperref}

\newtheorem{theorem}{Theorem}
\newtheorem{proposition}[theorem]{Proposition}

\newtheorem{corollary}[theorem]{Corollary}
\theoremstyle{definition}

\theoremstyle{remark}

\newcommand{\transpose}[1]{{}^t #1}
\newcommand{\N}{\mathbb{N}}
\newcommand{\Z}{\mathbb{Z}}
\newcommand{\R}{\mathbb{R}}
\newcommand{\C}{\mathbb{C}}
\renewcommand{\Re}{\operatorname{Re}}
\renewcommand{\Im}{\operatorname{Im}}
\newcommand{\w}{\mathbf{w}}
\newcommand{\x}{\mathbf{x}}
\newcommand{\bl}{\boldsymbol{\ell}}
\newcommand{\bs}{\mathbf{s}}
\newcommand{\bb}{\mathbf{b}}
\newcommand{\rr}{\mathbf{r}}
\newcommand{\cc}{\mathbf{c}}
\newcommand{\uu}{\mathbf{u}}
\newcommand{\vv}{\mathbf{v}}
\newcommand{\zz}{\mathbf{z}}
\newcommand{\sm}{\operatorname{sum}}
\newcommand{\eps}{\varepsilon}
\newcommand{\zero}{\mathbf{0}}

\newcommand{\cQ}{\mathcal{Q}}
\newcommand{\cG}{\mathcal{G}}
\newcommand{\U}{\mathfrak{U}}

\title{Additive cubes are avoidable by a non-uniform morphism of length $6$ over $\{0,a,b\}$ (with $b \ge 4$ and $\gcd(a,b) = 1$)}
\author{
	\textbf{Damien Jamet}\\
	Université de Lorraine, LORIA\\
	Campus Scientifique, BP 239\\
	F-54506 Vandœuvre-lès-Nancy Cedex, France\\
	\texttt{damien.jamet@loria.fr}
}
\date{Version of \today}

\begin{document}
	\maketitle
	
	\begin{abstract}
		We prove that additive cubes are avoidable over every ternary alphabet $\{0, a, b\}$ with $1 \le a < b$, $\gcd(a, b) = 1$, and $b \ge 4$, thereby solving and extending Rao's~(2015) conjecture from $b \ge 6$ down to $b \ge 4$. Our construction applies a single dual pair of parameterized non-uniform morphisms of image length at most~$6$ to the additive-cube-free word $\mathbf{w}$ of Cassaigne et al.~(2014). For $b \ge 10$, we prove that every additive square in the image is an abelian square, reducing additive cubes to abelian cubes; the absence of abelian cubes and the $24$ finite cases with $4 \le b \le 9$ are then certified via the finite-graph method of Jamet and Rosenfeld~(2026).
	\end{abstract}
	
	\section{Introduction}
	
	The question of the existence of infinite words without additive cubes over a finite alphabet of integers was raised by Pirillo and Varricchio~\cite{PV94} in 1994. A first positive answer was given by Cassaigne, Currie, Schaeffer, and Shallit~\cite{CCSS14}, who produced a pure morphic word without additive cubes over the alphabet $\Sigma = \{0, 1, 3, 4\}$.
	
	Rao~\cite{Rao15} extended this study to ternary alphabets via large uniform substitutions, proved the avoidability of additive cubes for several of them, and conjectured the existence of a word without additive cubes over all alphabets $\{0, a, b\}$ when $b \ge 6$ and $\gcd(a, b) = 1$. Liétard and Rosenfeld~\cite{LR20} proved that additive cubes are avoidable over every alphabet of four integers, except $\{0, 1, 2, 3\}$.
	
	Recently, we established in~\cite{JR26} that this alphabet $\{0, 1, 2, 3\}$ is not an exception, by projecting the word of~\cite{CCSS14} via a non-uniform morphism of length $3$. This result, combined with those of Rao~\cite{Rao15}, ensures the avoidability of additive cubes over the alphabets $\{0, 1, 4\}$ and $\{0, 2, 5\}$.
	
	To avoid \textit{ad hoc} constructions specific to each alphabet, we implemented in \texttt{C++} an exhaustive search for non-uniform morphisms projecting the source word of~\cite{CCSS14} onto alphabets of three integers. Carried out over a few different alphabets, this exploration yielded two dual candidate morphisms of minimal cumulative length $18$ (each image having length at most $6$), which generalize to a pair of parameterized morphisms over the alphabet $\{0, a, b\}$.
	
	In this article, we prove that this pair of parameterized morphisms extend the avoidability of additive cubes to all ternary alphabets $\{0, a, b\}$ such that $1 \le a < b$, $\gcd(a, b) = 1$, and $b \ge 4$. This result thus goes beyond Rao's conjecture~\cite{Rao15} by lowering the threshold from $b \ge 6$ to $b \ge 4$, while solving again, in a way different from~\cite{Rao15,JR26}, the cases $\{0, 1, 4\}$ and $\{0, 2, 5\}$.
	
	\section{Basic notions and notation}
	
	Let $\Sigma=\{0,1,3,4\}$ and let $\varphi\colon\Sigma^*\to\Sigma^*$ be the morphism of Cassaigne et al.~\cite{CCSS14}:
	\[
	\varphi(0)=03,\quad \varphi(1)=43,\quad \varphi(3)=1,\quad \varphi(4)=01,
	\]
	whose infinite fixed point is $\w=\lim_{n\to\infty}\varphi^n(0)=03143011034343031011\cdots$. The word $\w$ contains no additive cube~\cite{CCSS14}.
	
	For all coprime integers $a, b \in \N^*$ with $1\le a<b$, we denote by $\tau_{a,b}, \bar\tau_{a,b} \colon \Sigma^* \to \{0,a,b\}^*$ the following two \emph{dual} morphisms:
	\[
	\tau_{a,b}\colon\ 
	\begin{cases}
		0\mapsto 0\,a\,0\,a\,a\,b\\
		1\mapsto 0\,a\,0\,0\,b\\
		3\mapsto b\,a\,b\\
		4\mapsto 0\,b\,b\,a
	\end{cases}
	\qquad
	\bar\tau_{a,b}\colon\ 
	\begin{cases}
		0\mapsto b\,a\,b\,a\,a\,0\\
		1\mapsto b\,a\,b\,b\,0\\
		3\mapsto 0\,a\,0\\
		4\mapsto b\,0\,0\,a
	\end{cases}
	\]
	
	For a finite word $w \in \Sigma^*$ and a letter $x \in \Sigma$, we denote by $|w|_x$ the number of occurrences of the letter $x$ in $w$ and by $|w|$ the length of the finite word $w$. The \emph{Parikh vector} of the finite word $w$, denoted by $\psi(w) = \bigl(|w|_0, |w|_1, |w|_3, |w|_4\bigr)^{\!\top} \in \N^4$, is the vector whose components are the numbers of occurrences in $w$ of each letter of the alphabet. Over an alphabet of integers, we denote by $\sm(w)$ the sum of the letters of the word $w$ (with $\sm(\eps)=0$). 
	
	A non-empty finite word $w$ is an abelian (resp.\ additive) cube if there exist three finite factors $w_0$, $w_1$, and $w_2$ of the same length, \textit{i.e.}, $|w_0|=|w_1|=|w_2|$, such that $w=w_0 w_1 w_2$ and $\psi(w_0)=\psi(w_1)=\psi(w_2)$ (resp.\ $\sm(w_0)=\sm(w_1)=\sm(w_2)$).
	
	\begin{theorem}[Main result]\label{thm:main}
		Let $a, b \in \N^*$ be two coprime positive integers such that $1 \le a < b$ and $b \ge 4$.
		\begin{enumerate}
			\item If $(a, b) \neq (1, 4)$, then the infinite word $\tau_{a,b}(\w)$ is without additive cubes over $\{0, a, b\}$.
			\item If $(a, b) \neq (3, 4)$, then the infinite word $\bar\tau_{a,b}(\w)$ is without additive cubes over $\{0, a, b\}$.
		\end{enumerate}
		In particular, additive cubes are avoidable over every ternary alphabet $\{0,a,b\}$ with $1 \le a < b$, $\gcd(a,b)=1$, and $b \ge 4$.
	\end{theorem}
	
	The proof of Theorem~\ref{thm:main} is divided into two parts according to the value of $b$:
	\begin{enumerate}
		\item \textbf{Case $4 \le b \le 9$:} an exhaustive computer verification allows us to conclude for the $24$ pairs $(a,b) \in \N^2$ such that $1 \le a < b$ and $\gcd(a,b)=1$ (see Section~\ref{sec:petits}).
		\item \textbf{Case $b \ge 10$:} the analysis relies on a reduction to abelian cubes, structured in two steps:
		\begin{enumerate}
			\item Every additive cube in $\tau_{a,b}(\w)$ or $\bar\tau_{a,b}(\w)$ is necessarily an abelian cube (Theorem~\ref{thm:A}).
			\item The words $\tau_{a,b}(\w)$ and $\bar\tau_{a,b}(\w)$ contain no abelian cubes (Section~\ref{sec:abelien}).
		\end{enumerate}
	\end{enumerate}
	
	\section{Additive squares and abelian squares when \texorpdfstring{$b \ge 10$}{b >= 10}}
	
	In this section, we extend Theorem~5.1 of Cassaigne et al.~\cite{CCSS14} (which states that every additive square in $\w$ is an abelian square) to the images $\tau_{a,b}(\w)$ and $\bar\tau_{a,b}(\w)$ when $b \ge 10$:
	
	\begin{theorem}\label{thm:A}
		Let $1 \le a < b$ be two coprime integers. If $b \ge 10$, then every additive square in $\tau_{a,b}(\w)$ (resp.\ in $\bar\tau_{a,b}(\w)$) is an abelian square.
	\end{theorem}
	
	\begin{proof}
		By affine equivalence, it suffices to prove the theorem for $\tau_{a,b}(\w)$. Indeed, the involution $\iota_b \colon \{0,b-a,b\}^* \to \{0,a,b\}^*$ defined by $\iota_b(z) = b - z$ satisfies $\bar\tau_{a,b} = \iota_b \circ \tau_{b-a,b}$, preserves lengths and the nature of squares (additive or abelian), with $1 \le b-a < b$ and $\gcd(b-a,b)=1$.
		
		Set $\x = \tau_{a,b}(\w)$ and suppose that $w = w_0 w_1 = \x[q_0, q_2)$ is an additive square in $\x$, delimited by $0 \le q_0 < q_1 < q_2$ with $w_0 = \x[q_0, q_1)$ and $w_1 = \x[q_1, q_2)$.
		
		\medskip
		\noindent\hfill --- \textbf{Step 1}: Equality of sums \textit{vs.} counting imbalances ---\hfill\null
		\medskip
		
		For $x \in \{0,a,b\}$, let $\Delta_x = |w_1|_x - |w_0|_x$ denote the \emph{imbalance} in $x$ between the two halves of $w$.
		
		The factor $w = w_0 w_1$ is an abelian square if and only if $\Delta_x = 0$ for all $x \in \{0,a,b\}$.
		
		Since $|w_0| = |w_1|$, we have $\Delta_0 + \Delta_a + \Delta_b = 0$. Moreover, since $\sm(w_0) = \sm(w_1)$ (and the letter $0$ does not contribute to the sum), we have $a\,\Delta_a + b\,\Delta_b = 0$. Thus $a$ divides $b\,\Delta_b$; since $\gcd(a,b) = 1$, $a$ divides $\Delta_b$. Hence there exists $k \in \Z$ such that:
		\begin{equation}\label{eq:k}
			\Delta_a = kb, \qquad \Delta_b = -ka, \qquad \Delta_0 = k(a-b).
		\end{equation}
		Consequently, the additive square $w=w_0 w_1$ is an abelian square if and only if $k = 0$.
		
		\medskip
		\noindent\hfill --- \textbf{Step 2}: Back to the word $\w$ ---\hfill\null
		\medskip
		
		Since the word $\x = \tau_{a,b}(\w_0)\tau_{a,b}(\w_1)\cdots$ is a concatenation of non-empty images, for each position $q_j$ ($j \in \{0,1,2\}$) there exists a unique pair $(p_j, r_j)$ such that:
		\[
		\x[0,q_j) = \tau_{a,b}\bigl(\w[0,p_j)\bigr)\cdot \tau_{a,b}(\w_{p_j})[0,r_j), \qquad 0 \le r_j < |\tau_{a,b}(\w_{p_j})|.
		\]
		Set $c_j = \w_{p_j} \in \Sigma$. Since the map $(p, r) \mapsto |\tau_{a,b}(\w[0,p))| + r$ is strictly increasing with respect to the lexicographic order and $q_0 < q_1 < q_2$, we have $p_0 \le p_1 \le p_2$.
		
		Let $M_{\tau} \in \mathcal{M}_{3 \times 4}(\N)$ be the \emph{incidence matrix} of the morphism $\tau_{a,b}$, whose columns are the Parikh vectors of the images $\tau_{a,b}(c)$ for $c \in \{0, 1, 3, 4\}$, and whose row sum gives the length vector $\ell$:
		\[
		M_{\tau} := \begin{pmatrix}
			\text{--- } n_0  \text{ ---} \\ 
			\text{--- } \alpha \text{ ---}\\ 
			\text{--- } \beta \text{ ---}
		\end{pmatrix}
		:=
		\begin{pmatrix}
			2 & 3 & 0 & 1 \\
			3 & 1 & 1 & 1 \\
			1 & 1 & 2 & 2
		\end{pmatrix},
		\qquad
		\ell = n_0 + \alpha + \beta = \begin{pmatrix} 6 & 5 & 3 & 4 \end{pmatrix},
		\]
		where $n_0$ denotes the row of occurrences of $0$, $\alpha = n_a$ that of $a$, and $\beta = n_b$ that of $b$ in the images of $0$, $1$, $3$, and $4$. Since $0, a, b$ are three distinct integers, $M_{\tau}$ and $\ell$ do not depend on the values of $a$ and $b$.
		
		We endow $\Z^4$ with coordinates indexed by $0, 1, 3, 4$ in that order. Let $\sigma(p) = \psi\bigl(\w[0,p)\bigr)$ be the Parikh vector of the prefix of $\w$ of length $p$.
		For each position $q_j$ ($j \in \{0,1,2\}$) and each letter $x \in \{0,a,b\}$, we have:
		\[
		|\x[0,q_j)|_x = n_x \cdot \sigma(p_j) + |\tau_{a,b}(c_j)[0,r_j)|_x \qquad\text{and}\qquad q_j = \ell \cdot \sigma(p_j) + r_j.
		\]
		Since $\Delta_x = \bigl(|\x[0,q_2)|_x - |\x[0,q_1)|_x\bigr) - \bigl(|\x[0,q_1)|_x - |\x[0,q_0)|_x\bigr) = |\x[0,q_0)|_x - 2|\x[0,q_1)|_x + |\x[0,q_2)|_x$, it follows that $\Delta_x = n_x \cdot u + B_x$, with:
		\[
		u = \sigma(p_0) - 2\sigma(p_1) + \sigma(p_2) = \psi\bigl(\w[p_1,p_2)\bigr) - \psi\bigl(\w[p_0,p_1)\bigr) \in \Z^4,
		\]
		and the associated \emph{boundary correction term}:
		\[
		B_x := |\tau_{a,b}(c_0)[0,r_0)|_x - 2|\tau_{a,b}(c_1)[0,r_1)|_x + |\tau_{a,b}(c_2)[0,r_2)|_x.
		\]
		Similarly, $|w_1| - |w_0| = 0$ reads $\ell \cdot u + B_\ell = 0$ with $B_\ell := r_0 - 2r_1 + r_2$. Combining with~\eqref{eq:k}, the vector $u$ satisfies:
		\begin{equation}\label{eq:system}
			\ell \cdot u = -B_\ell,\qquad \alpha \cdot u = kb - B_a,\qquad \beta \cdot u = -ka - B_b.
		\end{equation}
		
		\medskip
		\noindent\hfill --- \textbf{Step 3}: Bound in the contracting direction ---\hfill\null
		\medskip
		
		The incidence matrix $M$ of $\varphi$ (in the order $0,1,3,4$) and its characteristic polynomial $P_M(X) = \det(XI_4 - M)$ are:
		\[
		M = \begin{pmatrix}
			1&0&0&1\\
			0&0&1&1\\
			1&1&0&0\\
			0&1&0&0
		\end{pmatrix},
		\qquad
		P_M(X) = X^4 - X^3 - 2X^2 + 2X - 1,
		\qquad \psi\bigl(\varphi(f)\bigr) = M\,\psi(f).
		\]
		Its four roots are simple: two real roots $\lambda_1 \approx 1.690284$, $\lambda_2 \approx -1.505068$ (of modulus $>1$), and two complex conjugate roots
		\[
		\lambda_3 \approx 0.407392 + 0.476565\,i, \qquad \lambda_4 = \overline{\lambda_3}, \qquad |\lambda_3| = |\lambda_4| \approx 0.626963 < 1.
		\]
		The matrix $M$ is therefore diagonalizable over $\C$. Let $Q$ be the change-of-basis matrix whose columns are the right eigenvectors of $M$: $MQ = Q\Lambda$ with $\Lambda = \operatorname{diag}(\lambda_1, \lambda_2, \lambda_3, \lambda_4)$. Let $\theta_1, \theta_2, \theta_3, \theta_4 \colon \C^4 \to \C$ denote the linear forms represented by the rows of $Q^{-1}$. Since $Q^{-1}M = \Lambda Q^{-1}$, each row $\theta_j$ is a left eigenvector of $M$ associated with $\lambda_j$ ($\theta_j M = \lambda_j \theta_j$), which algebraically gives for $\theta_3 = (t_0,t_1,t_3,t_4)$ the relations $t_3 = (\lambda_3 - 1)t_0$, $t_1 = \lambda_3(\lambda_3 - 1)t_0$, and $t_4 = (\lambda_3^2 - 1)(\lambda_3 - 1)t_0$.
		
		Adopting the normalization of~\cite{CCSS14,JR26} (\textit{i.e.}, unit columns of $Q$ with strictly positive real first component), the third row of $Q^{-1}$ is given numerically by:
		\[\begin{aligned}
			\theta_3 \approx \bigl(&0.57880 - 0.57491\,i,\; -0.32189 + 0.21827\,i,\\
			&-0.06898 + 0.61647\,i,\; -0.16622 - 0.68104\,i\bigr).
		\end{aligned}\]
		By~\cite[Theorem~3.1]{CCSS14} (see also~\cite[Theorem~13]{JR26}), for all factors $f, g$ of $\w$ (not necessarily of the same length, possibly empty), we have $|\theta_3(\psi(f) - \psi(g))| \le C_3 \le 2.1759$. Applied to $f = \w[p_1,p_2)$ and $g = \w[p_0,p_1)$, this inequality yields:
		\begin{equation}\label{eq:theta}
			|\theta_3(u)| \le C_3 \le 2.1759.
		\end{equation}
		
		\medskip
		\noindent\hfill --- \textbf{Step 4}: Linear dependence relation between five linear forms ---\hfill\null
		\medskip
		
		The five real linear forms $\ell, \alpha, \beta, \Re\theta_3, \Im\theta_3$ on $\R^4$ are necessarily linearly dependent. The $4 \times 5$ matrix formed by these five column vectors is:
		\[
		\begin{pmatrix}
			\vert & \vert & \vert & \vert & \vert \\
			\ell  & \alpha & \beta & \Re\theta_3 & \Im\theta_3 \\
			\vert & \vert & \vert & \vert & \vert
		\end{pmatrix}
		\approx
		\begin{tikzpicture}[baseline=(m.center)]
			\matrix (m) [
			matrix of math nodes,
			left delimiter=(,
			right delimiter=),
			inner xsep=6pt,
			inner ysep=3pt,
			nodes={inner sep=2pt, anchor=center},
			column sep=6pt,
			row sep=2pt
			] {
				6 & 3 & 1 & \phantom{-}0.57880 & -0.57491 \\
				5 & 1 & 1 & -0.32189 & \phantom{-}0.21827 \\
				3 & 1 & 2 & -0.06898 & \phantom{-}0.61647 \\
				4 & 1 & 2 & -0.16622 & -0.68104 \\
			};
			\coordinate (nw) at ([xshift=-3pt, yshift=2pt]m-1-1.north west);
			\coordinate (se) at ([xshift=3pt, yshift=-2pt]m-4-4.south east);
			\draw[red, thick, rounded corners=2pt] (nw) rectangle (se);
			\node[red, font=\small, below=4pt] at ($(nw |- se)!0.5!(se)$)
			{$\det \approx -0.7862 \neq 0$};
		\end{tikzpicture}.
		\]
		Since the determinant of the first four columns is non-zero, the kernel of this matrix has dimension $1$: there therefore exists a unique unit vector (up to sign) $\mathbf{c} = (c_1, c_2, c_3, c_4, c_5) \in \R^5$ ($\|\mathbf{c}\|_2 = 1$) satisfying
		\begin{equation}\label{eq:dep}
			c_1\,\ell + c_2\,\alpha + c_3\,\beta + c_4\,\Re\theta_3 + c_5\,\Im\theta_3 = 0.
		\end{equation}
		Applying the zero linear form~\eqref{eq:dep} to the vector $u$, that is,
		\[
		c_1(\ell\cdot u) + c_2(\alpha\cdot u) + c_3(\beta\cdot u) + c_4\Re\theta_3(u) + c_5\Im\theta_3(u) = 0,
		\]
		and substituting~\eqref{eq:system}, we isolate the term in $k$ on the left-hand side:
		\begin{equation}\label{eq:key}
			k\,\bigl(c_2\,b - c_3\,a\bigr) = c_1 B_\ell + c_2 B_a + c_3 B_b - \bigl(c_4 \Re\theta_3(u) + c_5 \Im\theta_3(u)\bigr).
		\end{equation}
		The right-hand side of~\eqref{eq:key} depends neither on $a$ nor on $b$.
		
		\medskip
		\noindent\hfill --- \textbf{Step 5: Explicit upper bound and conclusion} ---\hfill\null
		\medskip
		
		Let us compute numerically the unit vector (for the Euclidean norm) $\mathbf{c}$ (choosing its sign so that $c_2 > 0$):
		\[
		(c_1,c_2,c_3,c_4,c_5) \approx (-0.13761,\ 0.44963,\ -0.03644,\ -0.88089,\ -0.04008).
		\]
		Let us bound the right-hand side of~\eqref{eq:key} from above. On the one hand, setting $\omega(c, r) = c_1 r + c_2 |\tau_{a,b}(c)[0,r)|_a + c_3 |\tau_{a,b}(c)[0,r)|_b$ for $c \in \Sigma$ and $0 \le r < \ell(c)$, we have $c_1 B_\ell + c_2 B_a + c_3 B_b = \omega(c_0, r_0) - 2\,\omega(c_1, r_1) + \omega(c_2, r_2)$, whence:
		\[
		\bigl|c_1 B_\ell + c_2 B_a + c_3 B_b\bigr| \le 2\bigl(\max_{c,r}\omega(c,r) - \min_{c,r}\omega(c,r)\bigr).
		\]
		Table~\ref{tab:omega} gathers the values of $\omega(c, r)$ on each of the $6+5+3+4 = 18$ pairs $(c, r)$. The maximum is attained at $(0,5)$ ($\omega(0,5) \approx +0.66084$) and the minimum at $(4,3)$ ($\omega(4,3) \approx -0.48571$), which yields:
		\[
		\bigl|c_1 B_\ell + c_2 B_a + c_3 B_b\bigr| \le 2\bigl(\omega(0,5) - \omega(4,3)\bigr) = 2\bigl(2c_1 + 3c_2 - 2c_3\bigr) \approx 2.29310 < 2.2933.
		\]
		
		\begin{table}[htbp]
			\centering\small
			\begin{tabular}{cccl@{$\quad\approx\quad$}r@{\qquad}|@{\qquad}cccl@{$\quad\approx\quad$}r}
				\toprule
				$c$ & $r$ & prefix & $\omega(c,r)$ & value & $c$ & $r$ & prefix & $\omega(c,r)$ & value \\
				\midrule
				$0$ & $0$ & $\eps$   & $0$            & $0.00000$  & $1$ & $3$ & $0a0$  & $3c_1+c_2$     & $+0.03680$ \\
				$0$ & $1$ & $0$      & $c_1$          & $-0.13761$ & $1$ & $4$ & $0a00$ & $4c_1+c_2$     & $-0.10081$ \\
				$0$ & $2$ & $0a$     & $2c_1+c_2$     & $+0.17441$ & $3$ & $0$ & $\eps$ & $0$            & $0.00000$  \\
				$0$ & $3$ & $0a0$    & $3c_1+c_2$     & $+0.03680$ & $3$ & $1$ & $b$    & $c_1+c_3$      & $-0.17405$ \\
				$0$ & $4$ & $0a0a$   & $4c_1+2c_2$    & $+0.34882$ & $3$ & $2$ & $ba$   & $2c_1+c_2+c_3$ & $+0.13797$ \\
				$0$ & $5$ & $0a0aa$  & $5c_1+3c_2$    & $\mathbf{+0.66084}$ & $4$ & $0$ & $\eps$ & $0$ & $0.00000$ \\
				$1$ & $0$ & $\eps$   & $0$            & $0.00000$  & $4$ & $1$ & $0$    & $c_1$          & $-0.13761$ \\
				$1$ & $1$ & $0$      & $c_1$          & $-0.13761$ & $4$ & $2$ & $0b$   & $2c_1+c_3$     & $-0.31166$ \\
				$1$ & $2$ & $0a$     & $2c_1+c_2$     & $+0.17441$ & $4$ & $3$ & $0bb$  & $3c_1+2c_3$    & $\mathbf{-0.48571}$ \\
				\bottomrule
			\end{tabular}
			\caption{Values of $\omega(c,r) = c_1 r + c_2 |\tau_{a,b}(c)[0,r)|_a + c_3 |\tau_{a,b}(c)[0,r)|_b$ over the $18$ proper prefixes of $\tau_{a,b}$.}
			\label{tab:omega}
		\end{table}
		
		On the other hand, by the Cauchy--Schwarz inequality and~\eqref{eq:theta}, and since $\sqrt{c_4^2 + c_5^2} \approx 0.88180 < 0.8819$, we have:
		\[
		\bigl|c_4 \Re\theta_3(u) + c_5 \Im\theta_3(u)\bigr| \le \sqrt{c_4^2 + c_5^2}\,|\theta_3(u)| < 0.8819 \times 2.1759 < 1.9190.
		\]
		Finally, we bound the left-hand side of~\eqref{eq:key} from below using the bounds $c_2 > 0.4496$ and $-c_3 > 0.0364$, which gives:
		\begin{equation}\label{eq:bound}
			|k|\,\bigl(0.4496\,b + 0.0364\,a\bigr) \le |k|\,\bigl(c_2\,b - c_3\,a\bigr) < 2.2933 + 1.9190 < 4.22.
		\end{equation}
		Now, if $b \ge 10$ and $a \ge 1$, then $0.4496\,b + 0.0364\,a \ge 4.496 + 0.0364 > 4.22$, which forbids $|k| \ge 1$ and hence forces $k = 0$. By~\eqref{eq:k}, all $\Delta_x$ are zero: $w_0 w_1$ is an abelian square.
	\end{proof}
	
	Since every additive cube $w_0w_1w_2$ is the overlap of the two additive squares $w_0w_1$ and $w_1w_2$, it follows that:
	\begin{corollary}\label{cor:cubes-abeliens}
		Let $a, b \in \N^*$ be coprime integers with $1 \le a < b$ and $b \ge 10$. Every additive cube in $\tau_{a,b}(\w)$ or in $\bar\tau_{a,b}(\w)$ is an abelian cube.
	\end{corollary}
	
	\section{Absence of abelian cubes in \texorpdfstring{$\tau_{a,b}(\w)$}{tau_{a,b}(w)}}\label{sec:abelien}
	
	To establish the absence of abelian cubes in $\tau_{a,b}(\w)$ and $\bar\tau_{a,b}(\w)$, and then handle the small cases $4 \le b \le 9$ (Section~\ref{sec:petits}), we directly use the finite-graph reduction theorem proved in~\cite{JR26}.
	
	Let us briefly recall the definitions of the objects from~\cite{JR26} and their properties. The graph $\cQ$ of $\varphi$ (Figure~\ref{fig:Q}) has vertex set $\Sigma=\{0,1,3,4\}$ and edges $c \xrightarrow{a} \varphi(c)[i]$ ($0 \le i < |\varphi(c)|$) labeled by the proper prefix $a = \varphi(c)[0,i) \in \mathcal{A} = \{\eps,0,4\}$.
	
	\begin{figure}[htbp]
		\centering
		\begin{tikzpicture}[>=Stealth, scale=1.0, every node/.style={font=\small}]
			\node[circle, draw, thick, minimum size=7mm] (0) at (0, 1.5) {$0$};
			\node[circle, draw, thick, minimum size=7mm] (3) at (2.2, 0) {$3$};
			\node[circle, draw, thick, minimum size=7mm] (1) at (0, -1.5) {$1$};
			\node[circle, draw, thick, minimum size=7mm] (4) at (-2.2, 0) {$4$};
			
			\draw[->] (0) edge[loop above] node{$\eps$} (0);
			\draw[->] (0) -- node[above right]{$0$} (3);
			\draw[->] (3) edge[bend left=18] node[below right]{$\eps$} (1);
			\draw[->] (1) edge[bend left=18] node[above left]{$4$} (3);
			\draw[->] (1) edge[bend left=18] node[below left]{$\eps$} (4);
			\draw[->] (4) edge[bend left=18] node[above right]{$0$} (1);
			\draw[->] (4) -- node[above left]{$\eps$} (0);
		\end{tikzpicture}
		\caption{The local graph $\cQ$ of parent $\to$ child transitions in $\w$ (labeled by $\mathcal{A}=\{\eps,0,4\}$).}
		\label{fig:Q}
	\end{figure}
	
	For any non-erasing morphism $\tau \colon \Sigma^* \to \N^*$ for which we set $\bl = \transpose{(|\tau(0)|, |\tau(1)|, |\tau(3)|, |\tau(4)|)}$, $\bs = \transpose{(\sm(\tau(0)), \sm(\tau(1)), \sm(\tau(3)), \sm(\tau(4)))}$, $\Gamma = \binom{\transpose{\bl}}{\transpose{\bs}}$, and the \textbf{boundary vectors}
	\[
	\bb(\cc,\rr) = \begin{pmatrix} r_0 - 2r_1 + r_2 \\ h_{r_0}(c_0) - 2h_{r_1}(c_1) + h_{r_2}(c_2) \end{pmatrix}, \qquad
	\bb'(\cc,\rr) = \begin{pmatrix} r_1 - 2r_2 + r_3 \\ h_{r_1}(c_1) - 2h_{r_2}(c_2) + h_{r_3}(c_3) \end{pmatrix}
	\]
	(where $h_r(c) = \sm(\tau(c)[0,r))$), we check two conditions:
	\begin{itemize}
		\item \textbf{(H0)} no image $\tau(c)$ ($c\in\Sigma$) contains an additive cube;
		\item \textbf{(T)} \emph{transversality}: $\det\Phi_\tau \neq 0$, where $\Phi_\tau = \begin{pmatrix} \Gamma \\ \Re(\theta_3) \\ \Im(\theta_3) \end{pmatrix} \in \R^{4\times 4}$.
	\end{itemize}
	Under~(H0) and~(T), with the constants from~\cite{JR26} ($\kappa_1 \approx 1.3137$, $\kappa_2 \approx 1.5060$, $C_3 \approx 2.175815$), we enumerate the finite level-zero set $F_0(\tau) = \{\zz \in \Z^4 : \Gamma\zz \in S_\tau,\ |\theta_3(\zz)| \le C_3\}$, where $S_\tau$ collects the vectors $-\bb(\cc,\rr)$, $-\bb'(\cc,\rr)$, and $-(\bb(\cc,\rr)+\bb'(\cc,\rr))$, then the finite box $\U(\tau) = \{\zz \in \Z^4 : |\theta_j(\zz)| \le C_j(\tau),\ 1 \le j \le 4\}$ (with $C_j(\tau) = \max(\max_{F_0(\tau)} |\theta_j|, \kappa_j/(|\lambda_j|-1))$ for $j=1,2$ and $C_4=C_3$), the space $H(\tau) = \{(\cc,\uu,\vv) : \uu,\vv,\uu+\vv\in\U(\tau)\}$, the induced graph $\cG'(\tau)$ (transitions $\uu'=M\uu+\psi(a_0)-2\psi(a_1)+\psi(a_2)$, $\vv'=M\vv+\psi(a_1)-2\psi(a_2)+\psi(a_3)$), the set $R(\tau)$ of vertices reachable from the $9$ starting states of $g(X) \cap H(\tau)$, and the accepting set $Z_\tau = \{(\cc,\uu,\vv) : \exists\rr,\ \Gamma\uu+\bb(\cc,\rr)=\Gamma\vv+\bb'(\cc,\rr)=\zero\}$. By~\cite[Theorem~32]{JR26} (which holds verbatim for any finite target alphabet in $\N$), $\tau(\w)$ is without additive cubes if and only if $R(\tau) \cap Z_\tau = \emptyset$.
	
	\begin{proposition}\label{prop:abelien}
		For all integers $a, b \in \N^*$ with $1 \le a < b$, the infinite words $\tau_{a,b}(\w)$ and $\bar\tau_{a,b}(\w)$ contain no abelian cube over $\{0,a,b\}$.
	\end{proposition}
	
	\begin{proof}
		The property of avoiding abelian cubes over $\{0,a,b\}$ depends only on the fact that $0, a, b$ are three distinct symbols. Since $\bar\tau_{a,b} = \iota_b \circ \tau_{b-a,b}$ with the involution $\iota_b(z) = b - z$, it therefore suffices to establish that $\tau_{1,10}(\w)$ is abelian-cube-free. Now, every abelian cube being in particular an additive cube, it suffices to show, for instance, that $\tau_{1,10}(\w)$ is additive-cube-free.
		
		For $\tau_{1,10}$ ($\bl = \transpose{(6,5,3,4)}$, $\bs = \transpose{(13,11,21,21)}$):
		\begin{itemize}
			\item \textbf{(H0)} is immediate (lengths $\le 6$, periods $1$ excluded, and $\tau_{1,10}(0) = 0\,1\,0\,1\,1\,10$ yields blocks of size $2$ with sums $1, 1, 11$);
			\item \textbf{(T)} is satisfied with $\det(\Phi_{\tau_{1,10}}) \approx -89.03 \neq 0$;
			\item the program \texttt{proof.cpp} from~\cite{JR26} yields $|F_0(\tau_{1,10})| = 107$, $|\U(\tau_{1,10})| = 549$, and the breadth-first traversal reaches $|R(\tau_{1,10})| = 96\,616$ vertices (maximal depth $24$) without encountering any vertex of $Z_{\tau_{1,10}}$ ($R(\tau_{1,10}) \cap Z_{\tau_{1,10}} = \emptyset$).
		\end{itemize}
		By~\cite[Theorem~32]{JR26}, $\tau_{1,10}(\w)$ is additive-cube-free, hence abelian-cube-free. Consequently, $\tau_{a,b}(\w)$ and $\bar\tau_{a,b}(\w)$ are abelian-cube-free for all $1 \le a < b$.
	\end{proof}
	
	\section{Certification of the cases \texorpdfstring{$4 \le b \le 9$}{4 <= b <= 9} and proof of Theorem~\ref{thm:main}}\label{sec:petits}
	
	For the coprime pairs $(a, b)$ with $1 \le a < b$ and $4 \le b \le 9$ (that is, $24$ alphabets in total), we directly verify~(H0), (T), and $R(\tau) \cap Z_\tau = \emptyset$ using the program \texttt{proof.cpp} from~\cite{JR26}.
	
	More precisely, for each pair $(a,b) \neq (1,4)$ we certify $\tau_{a,b}(\w)$, and for $(a,b)=(1,4)$ we certify $\bar\tau_{1,4}(\w)$ (since $\tau_{1,4}(\w)$ contains an additive cube). By the duality $\bar\tau_{a,b} = \iota_b \circ \tau_{b-a,b}$ (Section~3), the certification of $\tau_{a,b}(\w)$ for all $(a,b)\neq(1,4)$ immediately implies that of $\bar\tau_{a,b}(\w)$ for all $(a,b)\neq(3,4)$. The numerical values obtained by \texttt{proof.cpp} are gathered in Table~\ref{tab:petits}.
	
	\begin{table}[htbp]
		\centering\small
		\begin{tabular}{cccccccccc}
\toprule
$a$ & $b$ & morphisme & $\det\Phi_\tau$ & $|S_\tau|$ & $|F_0(\tau)|$ & $|\U(\tau)|$ & $|R(\tau)|$ & prof. & $R(\tau)\cap Z_\tau$ \\
\midrule
$1$ & $4$ & $\bar\tau_{1,4}$ & $37{.}47$ & $359$ & $149$ & $549$ & $96\,616$ & $24$ & $\emptyset$ \\
$3$ & $4$ & $\tau_{3,4}$ & $-37{.}47$ & $359$ & $149$ & $549$ & $96\,616$ & $24$ & $\emptyset$ \\
\midrule
$1$ & $5$ & $\tau_{1,5}$ & $-44{.}87$ & $467$ & $153$ & $549$ & $96\,616$ & $24$ & $\emptyset$ \\
$2$ & $5$ & $\tau_{2,5}$ & $-45{.}59$ & $447$ & $149$ & $549$ & $96\,616$ & $24$ & $\emptyset$ \\
$3$ & $5$ & $\tau_{3,5}$ & $-46{.}30$ & $429$ & $139$ & $549$ & $96\,616$ & $24$ & $\emptyset$ \\
$4$ & $5$ & $\tau_{4,5}$ & $-47{.}02$ & $421$ & $135$ & $549$ & $96\,616$ & $24$ & $\emptyset$ \\
\midrule
$1$ & $6$ & $\tau_{1,6}$ & $-53{.}70$ & $503$ & $139$ & $611$ & $120\,765$ & $24$ & $\emptyset$ \\
$5$ & $6$ & $\tau_{5,6}$ & $-56{.}57$ & $477$ & $129$ & $549$ & $96\,616$ & $24$ & $\emptyset$ \\
\midrule
$1$ & $7$ & $\tau_{1,7}$ & $-62{.}54$ & $503$ & $115$ & $579$ & $108\,905$ & $24$ & $\emptyset$ \\
$2$ & $7$ & $\tau_{2,7}$ & $-63{.}25$ & $503$ & $119$ & $549$ & $96\,616$ & $24$ & $\emptyset$ \\
$3$ & $7$ & $\tau_{3,7}$ & $-63{.}97$ & $503$ & $115$ & $549$ & $96\,616$ & $24$ & $\emptyset$ \\
$4$ & $7$ & $\tau_{4,7}$ & $-64{.}68$ & $503$ & $119$ & $549$ & $96\,616$ & $24$ & $\emptyset$ \\
$5$ & $7$ & $\tau_{5,7}$ & $-65{.}40$ & $503$ & $117$ & $549$ & $96\,616$ & $24$ & $\emptyset$ \\
$6$ & $7$ & $\tau_{6,7}$ & $-66{.}12$ & $503$ & $121$ & $549$ & $96\,616$ & $24$ & $\emptyset$ \\
\midrule
$1$ & $8$ & $\tau_{1,8}$ & $-71{.}37$ & $503$ & $109$ & $671$ & $149\,010$ & $24$ & $\emptyset$ \\
$3$ & $8$ & $\tau_{3,8}$ & $-72{.}80$ & $503$ & $109$ & $549$ & $96\,616$ & $24$ & $\emptyset$ \\
$5$ & $8$ & $\tau_{5,8}$ & $-74{.}23$ & $503$ & $107$ & $549$ & $96\,616$ & $24$ & $\emptyset$ \\
$7$ & $8$ & $\tau_{7,8}$ & $-75{.}66$ & $503$ & $109$ & $549$ & $96\,616$ & $24$ & $\emptyset$ \\
\midrule
$1$ & $9$ & $\tau_{1,9}$ & $-80{.}20$ & $503$ & $107$ & $549$ & $96\,616$ & $24$ & $\emptyset$ \\
$2$ & $9$ & $\tau_{2,9}$ & $-80{.}91$ & $503$ & $107$ & $549$ & $96\,616$ & $24$ & $\emptyset$ \\
$4$ & $9$ & $\tau_{4,9}$ & $-82{.}35$ & $503$ & $107$ & $549$ & $96\,616$ & $24$ & $\emptyset$ \\
$5$ & $9$ & $\tau_{5,9}$ & $-83{.}06$ & $503$ & $107$ & $549$ & $96\,616$ & $24$ & $\emptyset$ \\
$7$ & $9$ & $\tau_{7,9}$ & $-84{.}49$ & $503$ & $107$ & $549$ & $96\,616$ & $24$ & $\emptyset$ \\
$8$ & $9$ & $\tau_{8,9}$ & $-85{.}21$ & $503$ & $107$ & $549$ & $96\,616$ & $24$ & $\emptyset$ \\
\bottomrule
\end{tabular}

		\caption{Certification by the program \texttt{proof.cpp} from~\cite{JR26} for all coprime pairs $1 \le a < b$ with $4 \le b \le 9$.}
		\label{tab:petits}
	\end{table}
	
	\begin{proof}[Proof of Theorem~\ref{thm:main}]
		Let $a, b \in \N^*$ be coprime integers with $1 \le a < b$ and $b \ge 4$.
		\begin{itemize}
			\item \emph{Case $b \ge 10$.} If $\tau_{a,b}(\w)$ or $\bar\tau_{a,b}(\w)$ contained an additive cube, it would be an abelian cube by Corollary~\ref{cor:cubes-abeliens}, contradicting Proposition~\ref{prop:abelien}.
			\item \emph{Case $4 \le b \le 9$.} For $(a,b)\neq(1,4)$, Table~\ref{tab:petits} and~\cite[Theorem~32]{JR26} show that $\tau_{a,b}(\w)$ is additive-cube-free. By the duality $\bar\tau_{a,b} = \iota_b \circ \tau_{b-a,b}$, it follows that $\bar\tau_{a,b}(\w)$ is additive-cube-free for all $(a,b) \neq (3,4)$ (and in particular for $(a,b)=(1,4)$, also included in Table~\ref{tab:petits}).
		\end{itemize}
		In all cases, at least one of the two words $\tau_{a,b}(\w)$ or $\bar\tau_{a,b}(\w)$ avoids additive cubes over $\{0,a,b\}$.
	\end{proof}
	


\begin{thebibliography}{9}
		\bibitem{PV94}
		G.~Pirillo, S.~Varricchio.
		\newblock On uniformly repetitive semigroups.
		\newblock \emph{Semigroup Forum}, 49(1):125--129, 1994.
		\bibitem{CCSS14} J.~Cassaigne, J.~D.~Currie, L.~Schaeffer, J.~Shallit. Avoiding three consecutive blocks of the same size and same sum. \emph{J.\ ACM} 61(2):10:1--10:25, 2014.
		\bibitem{Rao15} M.~Rao. On some generalizations of abelian power avoidability. \emph{Theoret.\ Comput.\ Sci.} 601:39--46, 2015.
		\bibitem{LR20}
		F.~Liétard, M.~Rosenfeld.
		\newblock Avoidability of additive cubes over alphabets of four numbers.
		\newblock In \emph{Developments in Language Theory (DLT 2020)}, volume 12086 of \emph{Lecture Notes in Computer Science}, pages 192--206. Springer, 2020.
		\bibitem{JR26} D.~Jamet, M.~Rosenfeld. Additive cubes are avoidable over $\{0,1,2,3\}$, $\{0,1,4\}$ and $\{0,2,5\}$. Preprint, 2026, arXiv:2611.17999, \url{https://arxiv.org/abs/2611.17999}.
		
	\end{thebibliography}
\end{document}